\documentclass[final,twoside,11pt]{entics}
\usepackage{enticsmacro}
\usepackage{graphicx}
\usepackage[all]{xy}
\usepackage{amsmath,amssymb}
\usepackage{tikz}
\usepackage{caption}
\newcommand{\ua}{\mathord{\uparrow}}
\newcommand{\da}{\mathord{\downarrow}}

\newcommand{\ddo}{\rotatebox[]{90}{$\twoheadleftarrow$}}
\newcommand{\uup}{\rotatebox[]{-90}{$\twoheadleftarrow$}}

\newcommand{\nll}{
\begin{tikzpicture}[baseline]
\draw (0,-0.2em) -- (1em,0.8em) ;
\node at (0.5em,0.3em) {$\ll$};
\end{tikzpicture}}

\def\conf{ISDT 2022} 	%%Fill in the acronym for your conference (with year)
\volume{2}			%Fill in the ENTICS volume number here
\def\volu{2}			% and here
\def\papno{4}			%Fill in your paper number here
\def\mydoi{{\normalshape  \href{https://doi.org/10.46298/entics.10344}{10.46298/entics.10344}}}
\def\copynames{H. Hou, Q. Li} %% Fill in the first initial and last name of the authors
\def\runtitle{Weakly meet $s_{Z}$-continuity and $\delta_{Z}$-continuity}

\def\lastname{Hou and Li}
\begin{document}
%%%Note the beginning and end of the frontmatter section that starts here%%%%%
\begin{frontmatter}
  \title{Weakly Meet $s_{Z}$-continuity and $\delta_{Z}$-continuity} %\tnoteref{t1}}
  %%%%Now the author(s) names(s)%%%%%
  \author{Huijun Hou\thanksref{ALL}\thanksref{myemail}}	%%Note NO SPACE between
   \author{Qingguo Li\thanksref{ALL}\thanksref{liqingguoli@aliyun.com}}		%last name and \thanksref{...}
    %%%Next come the addresses%%%%
   \address[a]{School of Mathematics\\ Hunan University\\				%or between \thanksrefs...
    Changsha, China}
    							
  \thanks[ALL]{This work is supported by the National Natural Science Foundation of China (No.12231007)}   %%Test of \thanks[ALL} here..
   \thanks[myemail]{Email: \href{houhuijun2021@163.com} {\texttt{houhuijun2021@163.com}}}
   %%%Note: if both authors share same institution, the only list the address once, after the second
   %%%author.
   %%%There also is a link from the first author to the co-author's address to show how to list
   %%%affiliations to more than one institution, when needed.
  \thanks[liqingguoli@aliyun.com]{Corresponding author, Email:  \href{mailto:liqingguoli@aliyun.com} {\texttt{\normalshape liqingguoli@aliyun.com}}}
\begin{abstract}
  Based on the concept of weakly meet $s_{Z}$-continuity put forward by Xu and Luo in \cite{qzm}, we further prove that if the subset system $Z$ satisfies certain conditions,  a poset is $s_{Z}$-continuous if and only if it is weakly meet $s_{Z}$-continuous and $s_{Z}$-quasicontinuous, which improves a related result   given by Ruan and Xu in \cite{sz}. Meanwhile, we provide a characterization for the poset to be weakly meet $s_{Z}$-continuous, that is, a poset with a lower hereditary $Z$-Scott topology is weakly meet $s_{Z}$-continuous if and only if it is locally weakly meet $s_{Z}$-continuous. In addition, we introduce a monad on the new category $\mathbf{POSET_{\delta}}$ and characterize its $Eilenberg$-$Moore$ algebras concretely.
\end{abstract}
\begin{keyword}
  Weakly meet $s_{Z}$-continuous poset \sep $\delta_{Z}$-continuous poset \sep Monad\sep $\mathrm{Eilenberg}$-$\mathrm{Moore}$ algebras
\end{keyword}

\end{frontmatter}
\section{Introduction}\label{intro}
Recall that the concept of subset system on the category $\mathbf{POSET}$ of posets was proposed by Wright et al. in \cite{ipic}. It originally  aimed at applying posets with $Z$-set structures to problems in computer science, particularly, to fixed point semantics for programming language. In addition, the set system includes many systems of sets which we are familiar with, such as directed sets, finite sets, connected sets and so on. Later, based on the suggestion given by Wright to study  the generalized counterpart of continuous poset by replacing directed sets with $Z$-sets, in \cite{zcp},  Baranga defined a kind of generalized way-below relation based on the $Z$-sets whose supremum exists. Furthermore, the author gave some characterizations for the $Z$-algebraic posets. Besides, Ern$\acute{e}$ introduced the concept of $s_{2}$-continuous posets by lending support to the cut operator of directed subsets instead of the existing sups, which is a pure order concept on posets, no longer depending on the dcpo.  Recently,  Zhang and Xu made use of the cut operator of directed sets again to define a new way-below relation between subsets, and then introduced $s_{2}$-quasicontinuous posets (see \cite{s2}). Combining with the notion of subset system, Xu and Luo in \cite{qzm} gave the definition of $s_{Z}$-quasicontinuous posets, and then, Ruan and Xu investigated  its properties in \cite{sz} concretely and mainly made the conclusion: when the subset system satisfies some conditions, a poset is $s_{Z}$-continuous if and only if it is $s_{Z}$-quasicontinuous and meet $s_{Z}$-continuous.

In this paper, we will first see that there is another characterization of $s_{Z}$-continuous posets. More precisely, in the case where the subset system $Z$ satisfies certain conditions,  a poset is $s_{Z}$-continuous if and only if it is $s_{Z}$-quasicontinuous and weakly meet $s_{Z}$-continuous, which is a  result stronger than that given by Xu in \cite{sz}, meanwhile, it reveals that the requirement of `$\sigma^{Z}(P) = \sigma_{Z}(P)$' is unnecessary. Then we focus on the weakly meet $s_{Z}$-continuity of a poset, and  show that a poset with a lower hereditary $Z$-Scott topology is weakly meet $s_{Z}$-continuous if and only if it is locally weakly meet $s_{Z}$-continuous. In order to investigate the $\Gamma$-faithful property, Ho and Zhao in \cite{zds} introduced a new beneath relation by the Scott closed subsets whose sups exist. Based on their work,  we find there exists a monad on $\mathbf{DCPO}$. Associating with the subset system, we will introduce a generalized beneath relation using the cuts of $Z$-Scott closed subsets, which is not necessary to consider whether the supremum exists. On this basis, we find a monad on the category $\mathbf{POSET_{\delta}}$ and characterize its Eilenberg-Moore algebras concretely.

\section{Preliminaries}

Let $P$ be a poset. For any $A\subseteq P,\ x\in P$, we write ${\uparrow} A = \{p\in P: p\geq$ a for some $a\in A\}$, ${\downarrow}$$A$ = $\{p\in P: p\leq a$ for some $a\in A\}$. In particular, $\uparrow$$x$ = $\uparrow$$\{x\}$ and $\downarrow$$x$ = $\downarrow$$\{x\}$. A subset $A\subseteq P$ is an \emph{upper set} (resp., a \emph{lower set}) if $A$ = $\uparrow$$A$ (resp., $A$ = $\downarrow$$A$). Let $A^{u}$ and $A^{l}$ denote the sets of all upper and lower bounds of $A$, respectively. The cut operator $\delta$ is defined by $E^{\delta} = E^{ul}$ for all $E\subseteq P$. Obviously, if the supremum of $E$ exists, then $x\in E^{\delta}$ iff $x\leq \sup E$. If $E\subseteq A$ for some subset $A$ of $P$, let $E^{\delta}\mid_{A} = \{p\in A: p\leq m\ \mathrm{for\ all}\ m\in E^{u}\cap A\}$. In particular, we write $E^{\delta}\mid_{\da m}$ as $E^{\delta}\mid_{m}$. We denote by $F\subseteq_{f} P$ if $F$ is a finite subset of poset $P$, and let $\mathbf{Fin} P = \{\ua F: F\subseteq_{f}P\}$. A mapping $\mathbf{min}: \mathbf{Fin} P\rightarrow 2^{P}$ is defined by $\mathbf{min} (\ua F) = \{x\in F: x\ \mathrm{is\ a\ minimal\ element\ of}\ F\}$.

For any $T_{0}$ space $X$, the partial order $\leq_{X}$ defined by $x\leq y$ iff $x$ is contained in the closure of $y$ is called the \emph{specialization order}. The topology on the poset $P$ generated by all principal filters $\ua x$ as a subbasis for the closed sets is called the $\emph{lower topology}$ and denoted by $\omega (P)$.

Let $\mathbf{POSET}$ denote the category of posets and monotone mappings. By \cite{ipic}, a $\emph{subset system}$  on $\mathbf{POSET}$ is a function $Z$ which assigns to each poset, a set  $Z(P)$ of subsets of $P$ such that
\begin{itemize}
 \item $\{x\}\in Z(P)$ for any $x\in P$, and
 \item if $f: P\rightarrow Q$ in $\mathbf{POSET}$ and $S\in Z(P)$, then $f(S)\in Z(Q)$.
\end{itemize}
$P$ is called a $\emph{Z-complete poset}$ ($zcpo$, for short), if $\sup D$ exists for each $D\in Z(P)$. A closure system on the set $X$ is a non-empty family $\mathcal{E}$ of subsets of $X$ which satisfies:
\begin{itemize}
\item $\bigcap_{i\in I}A_{i}\in \mathcal{E}$ for every nonempty family $\{A_{i}\}_{i\in I}\subseteq \mathcal{E}$, and
    \item $X\in \mathcal{E}$.

\end{itemize}

\begin{definition}
Let $P$ be a poset and let $\sigma ^{Z}(P) = \{U\subseteq P: \mathrm{for\ all}\ S\in Z(P), S^{\delta}\cap U\neq \emptyset \Rightarrow S\cap U\neq \emptyset\}$. The topology generated by the subbasic open subsets $\sigma^{Z}(P)$ is called $Z$-Scott topology on $P$ and denoted by $\sigma_{Z}(P)$.
\end{definition}

Let $\Gamma^{Z}(P) = \{A\subseteq P: \mathrm{for\ all}\ S\in Z(P), S\subseteq P\Rightarrow S^{\delta}\subseteq P\}$, obviously, $\Gamma^{Z}(P)$ is a subbasis for the closed subsets with respect to $Z$-Scott topology. We use $\Gamma_{Z}(P)$ to denote the set composed of all closed subsets regarding $Z$-Scott topology. Note that for any $U\in \sigma ^{Z}(P)(A\in \Gamma^{Z}(P)), U  = \ua U (A = \da A)$, so the definition above is the same as that given in \cite{sz}. Besides, the family $\Gamma^{Z}(P)$ and $\Gamma_{Z}(P)$ are closure systems on $P$, and the closure operators can be defined as follows:  For any $M\subseteq P$, $cl_{\sigma^{Z}(P)}(M) = \bigcap \{A\in \Gamma^{Z}(P): M\subseteq A\}$, $cl_{\sigma_{Z}(P)}(M) = \bigcap\{B\in \Gamma_{Z}(P): M\subseteq B\}$.

\begin{definition}(\cite{sz})
Let $P$ be a poset and $x\in P,\ A, B\subseteq P$.
\begin{enumerate}
\item $A$ is called $Z$-way below $B$, denoted by $A\ll_{Z} B$, if for any $S\in Z(P), S^{\delta}\cap \ua B\neq \emptyset$ implies $S\cap \ua A\neq \emptyset$. $F\ll_{Z}\{x\}$ is shortly written as $F\ll_{Z}x$. Let $\omega_{Z}(x) = \{F\subseteq_{f}P: F\ll_{Z}x\}, \Uparrow_{Z}A = \{x\in P: A\ll_{Z}x\},\ \uup_{Z}A = \{p\in P: a\ll_{Z}p\ \mathrm{for\ some}\ a\in A\},\ \ddo_{Z} x = \{y\in P: y\ll_{Z}x\}$. Specifically, we write $\ddo_{Z}^{x}y = \{m\in \da x: m\ll_{Z}y\ \mathrm{in}\ \da x\}$.
\item $P$ is called a weak $s_{Z}$-continuous poset, if for all $x\in P, x\in (\ddo_{Z}x)^{\delta}$. In addition, if $\ddo_{Z}x\in I_{Z}(P) = \{\da S: S\in Z(P)\}$, then $P$ is called $s_{Z}$-continuous.
\item $P$ is called an $s_{Z}$-quasicontinuous poset, if for all $p\in P$, $\{\ua F: F\in \omega_{Z}(P)\}\in Z(\mathbf{Fin} P)$ and $\ua p = \bigcap\{\ua F: F\in \omega_{Z}(P)\}$.
\end{enumerate}
\end{definition}

\begin{definition}(\cite{cscp})
A subset system $Z$ is hereditary if for any order embedding $f: P\rightarrow Q$ (that is, for any $x, y\in P, f(x)\leq f(y) \Leftrightarrow x\leq y$), $D\subseteq P$, $D\in Z(P)$ if and only if $f(D)\in Z(Q)$.
\end{definition}

\begin{definition}(\cite{sz})
Let $Z$ be a subset system.
\begin{enumerate}
\item $Z$ is called $union\ complete$, if for any poset $P$, $\mathcal{S}\in Z(Z(P))$, we have $\bigcup \mathcal{S}\in Z(P)$.
\item $Z$ is said to have the $finite\ family\ union\ property$, if for any poset $P$, $\{\mathcal{S}_{1}, \mathcal{S}_{2}, ...,\mathcal{S}_{n}\}\subseteq_{f} Z(\mathbf{Fin} P)$, we have $\{\bigcup_{i = 1}^{n}A_{i}: A_{i}\in \mathcal{S}_{i}, i = 1,2,...,n\}\in Z(\mathbf{Fin} P)$.
\item $Z$ ia said to have the property $M$, if for any poset $P$, $\ua F\in \mathbf{Fin} P$, we have $\da_{\mathbf{Fin} P}\ua F = \{\ua G\in \mathbf{Fin} P: \ua F\subseteq \ua G\}\in Z(\mathbf{Fin} P)$.
\end{enumerate}
\end{definition}

\begin{definition}(\cite{sz})
A subset system $Z$ is said to have the $Rudin\ property$, if for any poset $P$, $E = \ua E\subseteq P$, $\mathcal{G}\in Z(\mathbf{Fin} P)$, $\emptyset\notin \mathcal{G}$, and $\bigcap \mathcal{G}\subseteq E$. Then there exists $K\subseteq \bigcup\{\mathbf{min}(G): G\in \mathcal{G}\}$ such that
\begin{enumerate}
  \item[(i)] for any $G\in \mathcal{G}$, $K\cap \mathbf{min}(G)\neq\emptyset$,
  \item[(ii)] $K\in Z(P)$,
  \item[(iii)] $\bigcap\{\ua k: k\in K\}\subseteq E$, and
  \item[(iv)] for any $G, H\in \mathcal{G}$, $G\subseteq H$ implies $K\cap \mathbf{min} (G)\subseteq \ua (K\cap \mathbf{min} (H))$.
\end{enumerate}
$Z$ is called a $Rudin\ subset\ system$, if $Z$ is union-complete and possesses the Rudin property.

%\begin{enumerate}

%\end{enumerate}
\end{definition}

\section{Weakly meet $s_{Z}$-continuous posets}
\begin{definition}(\cite{sz})
$P$ is called weakly meet $s_{Z}$-continuous if for all $x\in P$ and all $D\in Z(P)$ with $x\in D^{\delta}$, we have $x\in cl_{\sigma^{Z}(P)}(\da x\cap \da D)$; $P$ is called meet $s_{Z}$-continuous if for all $x\in P$ and all $D\in Z(P)$ with $x\in D^{\delta}$, we have $x\in cl_{\sigma_{Z}(P)}(\da x\cap \da D)$.
\end{definition}

\begin{lemma}\label{wmc}
Let $P$ be a poset. The following conditions are equivalent:
\begin{enumerate}
\item[$\mathrm{(1)}$] $P$ is weakly meet $s_{Z}$-continuous.
\item[$\mathrm{(2)}$] For any $x\in P$ and any $U\in \sigma^{Z}(P)$, $\ua (\da x\cap U)\in \sigma^{Z}(P)$.
\end{enumerate}
\end{lemma}
\begin{proof}
$(1)\Rightarrow (2)$ Assume that $D\in Z(P)$, and $D^{\delta}\cap \ua (\da x\cap U)\neq \emptyset$. Then there exists an $m\in D^{\delta}$ with $m\in U$ and $m\leq x$. Since $P$ is weakly meet $s_{Z}$-continuous, we have $m\in cl_{\sigma^{Z}(P)}(\da m\cap \da D)$, which implies that $\da D\cap \da m\cap U\neq \emptyset$. Thus $\da D\cap \da x\cap U\neq \emptyset$ by $m\leq x$. So $D\cap \ua(\da x\cap U)\neq \emptyset$ and $\ua(\da x\cap U)\in \sigma^{Z}(P)$ holds.

$(2)\Rightarrow (1)$ For any $x\in P$, $D\in Z(P)$, if $x\in D^{\delta}$ and there is a $U\in \sigma^{Z}(P)$ such that $x\in U$, then by (2), $\ua(\da x\cap U)\in \sigma^{Z}(P)$. Since $x\in D^{\delta}\cap \ua(\da x\cap U)\neq \emptyset$, we have $D\cap \ua(\da x\cap U)\neq \emptyset$, this means $\da x\cap U\cap \da D\neq\emptyset$. So $x\in cl_{\sigma^{Z}(P)}(\da x\cap\da D)$.
\end{proof}

\begin{lemma}\label{2}
Let $P$ be a $Z$-complete semilattice. The following conditions are equivalent:
\begin{enumerate}
\item[$\mathrm{(1)}$] $P$ is weakly meet $s_{Z}$-continuous;
\item[$\mathrm{(2)}$] For any $x\in P$, $D\in Z(P)$, $x\wedge\vee D = \vee\{x\wedge d: d\in D\}$.
\end{enumerate}
\end{lemma}
\begin{proof}
$(1)\Rightarrow (2)$ We first claim that $y = \vee(\da y\cap \da D)$ if $y\in D^{\delta}$. It is obvious that $y$ is an upper bound of $\da y\cap \da D$. Suppose $z$ is also an upper bound of $\da y\cap \da D$ and $y\nleqslant z$, that is, $y\in P\setminus \da z$. Since $y\in cl_{\sigma^{Z}(P)}(\da y\cap \da D)$ by (1) and $P\setminus \da z\in \sigma^{Z}(P)$, we have $(P\setminus \da z)\cap\da y\cap \da D\neq\emptyset$. But this contracts  the fact that $\da y\cap \da D\subseteq \da z$. Thus $y = \vee(\da y\cap \da D)$. Now let $y_{0} = x\wedge\vee D$, then $y_{0}\in D^{\delta}$, which implies $y_{0} = \vee(\da y_{0}\cap \da D)$. Since $\da y_{0}\cap \da D = \da \{x\wedge d: d\in D\}$, we have $y_{0} = \vee\{x\wedge d: d\in D\}$, that is, $x\wedge\vee D = \vee\{x\wedge d: d\in D\}$.

$(2)\Rightarrow (1)$ For any $x\in P$, $U\in \sigma^{Z}(P)$, we need to prove $\ua(\da x\cap U)\in \sigma^{Z}(P)$. Assume  $D\in Z(P)$ with $D^{\delta}\cap \ua(\da x\cap U)\neq\emptyset$. Then there exists an $m\in U$, $m\leq x$ and $m\in D^{\delta}$. Thus $m\leq\vee D$ and $m = m\wedge\vee D = \vee\{m\wedge d: d\in D\}\in U$ by (2). Now for $m$, we define a monotone mapping $\varphi: P\rightarrow P$ by $\varphi(p) = m\wedge p$. Then $\varphi(D) = \{m\wedge d: d\in D\}\in Z(P)$. Hence, $m\wedge d_{0}\in U$ for some $d_{0}\in D$ as $U\in \sigma^{Z}(P)$, which implies that $D\cap \ua(\da x\cap U)\neq\emptyset$, that is, $\ua(\da x\cap U)\in \sigma^{Z}(P)$. So $P$ is weakly meet $s_{Z}$-continuous by Lemma \ref{wmc}.
\end{proof}

\begin{proposition}
Let $P$ be a poset. The following conditions are equivalent:
\begin{enumerate}
\item[$\mathrm{(1)}$] $P$ is weakly meet $s_{Z}$-continuous;
\item[$\mathrm{(2)}$] $\Gamma^{Z}(P)$ is weakly meet $s_{Z}$-continuous.
\end{enumerate}
\end{proposition}
\begin{proof}
$(1)\Rightarrow (2)$ By Lemma \ref{2}, we only need to prove that for any $A\in \Gamma^{Z}(P), \mathcal{D}\in Z(\Gamma^{Z}(P))$, $A\wedge (\vee \mathcal{D}) = \vee\{A\wedge D: D\in \mathcal{D}\}$, that is, $A\cap cl_{\sigma^{Z}(P)}(\bigcup \mathcal{D}) = cl_{\sigma^{Z}(P)}(\bigcup\{A\cap D: D\in \mathcal{D}\})$. Assume $x\in A\cap cl_{\sigma^{Z}(P)}(\bigcup \mathcal{D})$ and $U\in \sigma^{Z}(P)$ with $x\in U$. Then we have $\ua(\da x\cap U)\in \sigma^{Z}(P)$ by Lemma \ref{wmc}. As $x\in cl_{\sigma^{Z}(P)}(\bigcup \mathcal{D})$, $\ua(\da x\cap U)\cap D_{0}\neq\emptyset$ for some $D_{0}\in \mathcal{D}$, this means $A\cap U\cap D_{0}\neq\emptyset$ since $x\in A$ and $A$ is a lower set. Moreover, $(\bigcup\{A\cap D: D\in \mathcal{D}\})\cap U\neq\emptyset$. So $x\in cl_{\sigma^{Z}(P)}(\cup\{A\cap D: D\in \mathcal{D}\})$, and $A\cap cl_{\sigma^{Z}(P)}(\bigcup \mathcal{D}) \subseteq cl_{\sigma^{Z}(P)}(\bigcup\{A\cap D: D\in \mathcal{D}\})$ holds. Obviously, the conversely inclusion holds.

$(2)\Rightarrow (1)$ It is sufficient to prove that $\ua(\da x\cap U)\in \sigma^{Z}(P)$ for any $x\in P$, $U\in \sigma^{Z}(P)$. Let $D\in Z(P)$ with $D^{\delta}\cap \ua(\da x\cap U)\neq\emptyset$. Then there exists an $m\in \da x\cap U$ such that $m\in D^{\delta}$, which implies $\da m \in \{\da d: d\in D\}^{\delta}$. In addition, we know $\{\da d: d\in D\}\in Z(\Gamma^{Z}(P))$ since the mapping $\psi: P\rightarrow \Gamma^{Z}(P)$ defined by $\psi(p) = \da p$ is monotone. As $\Gamma^{Z}(P)$ is weakly meet $s_{Z}$-continuous, we have $\da m\in cl_{\sigma^{Z}(\Gamma^{Z}(P))}(\da \{\da m\}\cap \da\{\da d: d\in D\} )$. It is easy to verify that $\lozenge U = \{A\in \Gamma^{Z}(P): A\cap U\neq\emptyset\}\in \sigma^{Z}(\Gamma^{Z}(P))$ and $\da m\in \lozenge U$. So $\lozenge U\cap \da \{\da m\}\cap \da\{\da d: d\in D\} \neq\emptyset$, that is,  $C\in \Gamma^{Z}(P)$ belongs to this intersection. Moreover, there exists an element $c\in C\cap U$ satisfying $c\leq m\leq x$ and $c\leq d_{0}$ for some $d_{0}\in D$, this means $D\cap \ua(\da x\cap U)\neq\emptyset$. Hence, $\ua(\da x\cap U)\in \sigma^{Z}(P)$, and $P$ is weakly meet $s_{Z}$-continuous.
\end{proof}

\begin{lemma}\label{3}
Let $P$ be a weakly meet $s_{Z}$-continuous poset. If $F$ is a finite subset of $P$, then $int_{\sigma^{Z}(P)}(\ua F)\subseteq\cup\{\uup_{Z}x: x\in F\}$.
\end{lemma}
\begin{proof}
Suppose $F = \{x_{1}, x_{2}, ... , x_{n}\}$ and there exists an element $a\in int_{\sigma^{Z}(P)}(\ua F)$, but $a\notin \cup\{\uup x_{i}: i = 1,2,...,n\}$. Then $x_{i}\nll_{Z} a$ for any $x_{i}\in F$, that is, there exists $D_{i}\in Z(P)$ such that $a\in D_{i}^{\delta}$, but $x_{i}\notin \da D_{i}$, for $i = 1,2,...,n$. For $D_{1}\in Z(P)$ with $a\in D_{1}^{\delta}$, $a\in cl_{\sigma^{Z}(P)}(\da a\cap \da D_{1})$ by weakly meet $s_{Z}$-continuity. Then $int_{\sigma^{Z}(P)}(\ua F)\cap \da a\cap \da D_{1}\neq\emptyset$, which implies that there is a $y_{1}\in int_{\sigma^{Z}(P)}(\ua F)\cap \da a\cap \da D_{1}$. By $y_{1}\leq a$ and $a\in D_{2}^{\delta}$, we have $y_{1}\in D_{2}^{\delta}$. Similarly, we get that $y_{1}\in cl_{\sigma^{Z}(P)}(\da y_{1}\cap \da D_{2})$ and $int_{\sigma^{Z}(P)}(\ua F)\cap\da y_{1}\cap \da D_{2}\neq\emptyset$. So there is a $y_{2}\in int_{\sigma^{Z}(P)}(\ua F)\cap\da y_{1}\cap \da D_{2}$. By  induction, we find $y_{n}\in int_{\sigma^{Z}(P)}(\ua F)\cap\da y_{n-1}\cap \da D_{n}$, where $y_{0} = a$,  clearly, $y_{n}\in \bigcap_{i = 1}^{i = n}\da D_{i}$. Since $y_{n}\in int_{\sigma^{Z}(P)}(\ua F)\subseteq \ua F$, $y_{n}\geq x_{i_{0}}$ for some $i_{0}\in \{1,2,...,n\}$, this implies $x_{i_{0}}\in \da D_{i_{0}}$, which contradicts  that $x_{i_{0}}\notin \da D_{i_{0}}$. Hence, $int_{\sigma^{Z}(P)}(\ua F)\subseteq\cup\{\uup_{Z}x: x\in F\}$.
\end{proof}

\begin{lemma}(\cite{sz})\label{4}
Let $Z$ be a Rudin subset system which has the finite family union property and $P$ an $s_{Z}$-quasicontinuous poset. Then the following statements hold.
\begin{enumerate}
\item[$\mathrm{(1)}$] For any finite set $F$ in $P$, $\Uparrow_{Z}F\in \sigma^{Z}(P)$.
\item[$\mathrm{(2)}$] If $U\subseteq P$, then $U\in \sigma^{Z}(P)$ if and only if for any $x\in U$, there exists $F\subseteq_{f} P$ such that $x\in \Uparrow_{Z}F\subseteq\ua F\subseteq U$.
\end{enumerate}
\end{lemma}
\begin{lemma}
Let $Z$ be a Rudin subset system which has the finite family union property. If $P$ is weakly meet $s_{Z}$-continuous and $s_{Z}$-quasicontinuous, then for any finite subset $F$ of $P$, we have
\begin{center}
$\Uparrow_{Z}F = \uup_{Z} F$.
\end{center}
\end{lemma}
\begin{proof}
By Lemma \ref{3} and Lemma \ref{4}, obviously, $\Uparrow_{Z}F \subseteq \uup_{Z} F$. And the reverse containment is easy to verify, so we omit the proof.
\end{proof}

\begin{proposition}(\cite{qzm})\label{p1}
Let $Z$ be a Rudin subset system which possesses $M$ property. If $P$ is an $s_{Z}$-continuous poset, then $P$ is $s_{Z}$-quasicontinuous, and for any $p\in P$, $\omega_{Z}(p) = \{F\subseteq_{f}P: \exists y\ll_{Z}p, \mathrm{such\ that}\ y\in \ua F\}$.
\end{proposition}

\begin{proposition}(\cite{sz})\label{p2}
Let $P$ be an $s_{Z}$-continuous poset. Then $P$ is weakly meet $s_{Z}$-continuous.
\end{proposition}

\begin{proposition}\label{p3}
Let $P$ be a weakly meet $s_{Z}$-continuous poset. If for any $x,y\in P$, $x\nleq y$, there are $U\in \sigma^{Z}(P)$, $V\in \omega(P)$ such that $x\in U$, $y\in V$ and $U\cap V =\emptyset$, then $P$ is weak $s_{Z}$-continuous.
\end{proposition}
\begin{proof}
It suffices to prove that $x\in (\ddo_{Z}x)^{\delta}$ for any $x\in P$. Suppose that there is a $y\in (\ddo_{z}x)^{u}$ but $x\nleq y$. Then there are $U\in \sigma^{Z}(P)$, $V = P\setminus\ua F\in \omega(P)$ such that $x\in U$, $y\in V$ and $U\cap V =\emptyset$, so $U\subseteq\ua F$. Since $\ua(\da x\cap U)\in \sigma^{Z}(P)$ by Lemma \ref{wmc} and $x\in \ua(\da x\cap U)\subseteq\ua F$, we have $x\in int_{\sigma^{Z}(P)}(\ua F)\subseteq \uup_{Z}F$. Thus there is an $m\in F$ such that $m\in \ddo_{Z}x$. It follows that $m\leq y$, then $y\in \ua F$. But this contradicts that $y\in V$.
\end{proof}

\begin{theorem}
Let P be a poset and $Z$ a Rudin subset system which possesses the finite family union property and $M$ property. If $\ddo_{Z}x\in I_{Z}(P)$ for each $x\in P$, then the following conditions are equivalent:
\begin{enumerate}
\item[$\mathrm{(1)}$] $P$ is $s_{Z}$-continuous;
\item[$\mathrm{(2)}$] $P$ is weakly meet $s_{Z}$-continuous and $s_{Z}$-quasicontinuous;
\item[$\mathrm{(3)}$] $P$ is weakly meet $s_{Z}$-continuous, and for any $x\nleq y$ in $P$, there are $U\in \sigma^{Z}(P)$, $V\in \omega(P)$ such that $x\in U$, $y\in V$ and $U\cap V= \emptyset$.
\end{enumerate}
\end{theorem}
\begin{proof}
$(1)\Rightarrow (2)$ Straightforward by Proposition \ref{p1} and Proposition \ref{p2}.

$(2)\Rightarrow (3)$ For any $x\nleq y$, that is, $y\notin \ua x$, there is an $F\in \omega_{Z}(x)$ such that $y\notin \ua F$ by (2). So we get that there are $\Uparrow_{Z}F\in \sigma^{Z}(P)$, $P\setminus\ua F\in \omega(P)$ containing $x$ and $y$, respectively, and $\Uparrow_{Z}F\cap P\setminus\ua F = \emptyset$.

$(3)\Rightarrow (1)$ By Proposition \ref{p3} and $\ddo_{Z}x\in I_{Z}(P)$, we know $P$ is $s_{Z}$-continuous.
\end{proof}

\section{Posets with lower hereditary Z-Scott topology}

\begin{definition}
Let $P$ be a poset. The $Z$-Scott topology on $P$ is called $lower\ hereditary$ if for each closed subbasis $A$ of $P$, the $Z$-Scott topology of poset $A$ is precisely generated by the subbasic closed subsets of the form $B\cap A$, where $B\in \Gamma^{Z}(P)$, that is, $\Gamma^{Z}(A) = \{B\cap A: B\in \Gamma^{Z}(P)\}$.
\end{definition}

\begin{definition}
Let $P$, $Q$ be two posets. A mapping $f: P\rightarrow Q$ is called $\sigma^{Z}$-continuous if for any $A\in \Gamma^{Z}(Q)$, $f^{-1}(A)\in \Gamma^{Z}(P)$.

It is obvious that $f$ is monotone if $f$ is $\sigma^{Z}$-continuous.
\end{definition}

\begin{lemma}
Let $P$ and $Q$ be two posets and $f: P\rightarrow Q$. Consider the following three conditions:
\begin{enumerate}
\item[$\mathrm{(1)}$] $f$ is $\sigma^{Z}$-continuous.
\item[$\mathrm{(2)}$] For any $D\in Z(P)$, $f(D^{\delta}) \subseteq f(D)^{\delta}$.
\item[$\mathrm{(3)}$] $f(cl_{\sigma^{Z}(P)}(A)) \subseteq cl_{\sigma^{Z}(P)}(f(A))$ for each $A\subseteq P$.
\end{enumerate}
Then $(1)\Leftrightarrow(2)\Rightarrow(3)$.
\end{lemma}
\begin{proof}
Straightforward.
\end{proof}

\begin{lemma}\label{lh}
Let $P$ be a poset and $Z$ a subset hereditary subset system. Consider the following conditions:
\begin{enumerate}
\item[$\mathrm{(1)}$] The $Z$-Scott topology on $P$ is lower hereditary.
\item[$\mathrm{(2)}$] The inclusion map $i: \da x\rightarrow P$ is $\sigma^{Z}$-continuous for any $x\in P$.
\item[$\mathrm{(3)}$] For any $x\in P$ and $D\in Z(\da x)$, $D^{\delta}\mid_{x} = D^{\delta}$.
\item[$\mathrm{(4)}$] For any $A\in \Gamma^{Z}(P)$ and $D\in Z(A)$, $D^{\delta}\mid_{A} = D^{\delta}$.
\item[$\mathrm{(5)}$] For any $D\in Z(P)$, $D^{u}$ is filtered.
\end{enumerate}
Then $(5)\Rightarrow(1)\Leftrightarrow(2)\Leftrightarrow(3)\Leftrightarrow(4)$.
\end{lemma}
\begin{proof}
It is easy to verify that $(1)\Rightarrow(2)\Rightarrow(3)$.

$(3)\Rightarrow(4):$ It is clear that $D^{\delta}\subseteq D^{\delta}\mid_{A}$. Assume $m\in D^{u}\mid_{A}$. Then $D\subseteq \da m$ and $m\in A$. Since $Z$ is subset hereditary, $D\in Z(\da m)$. Thus we have $D^{\delta}\mid_{m} = D^{\delta}$ by $(3)$. Now we only need to prove that $D^{\delta}\mid_{A}\subseteq D^{\delta}\mid_{m}$. Assume $a\in D^{\delta}\mid_{A}$, $b\in D^{u}\mid_{m}$. Then $b\leq m$ and $b\in A$ as $m\in A$, which implies that $b\in D^{u}\mid_{A}$, so $a\leq b$. Hence, $D^{\delta}\mid_{A}\subseteq D^{\delta}\mid_{m}$.

$(4)\Rightarrow(1):$ We want to prove that $\Gamma^{Z}(A) = \{A\cap C: C\in \Gamma^{Z}(P)\}$ for any $A\in \Gamma^{Z}(P)$. For each $B\in \Gamma^{Z}(A)$, let $D\in Z(P)$ and $D\subseteq B$. Then $D\in Z(A)$ because $Z$ is subset hereditary. It follows that $D^{\delta}\mid_{A}\subseteq B$, which means $D^{\delta}\subseteq B$ since $D^{\delta}\mid_{A} = D^{\delta}$. Thus $B\in \Gamma^{Z}(P)$ and $\Gamma^{Z}(A)\subseteq \{A\cap C: C\in \Gamma^{Z}(P)\}$. Conversely, for any $C\in \Gamma^{Z}(P)$, let $D\in Z(A)$ with $D\subseteq A\cap C$. Then $D\in Z(P)$ and $D^{\delta}\subseteq A\cap C$ since $A\cap C\in\Gamma^{Z}(P)$. This implies that  $D^{\delta}\mid_{A}\subseteq A\cap C$. So $A\cap C\in \Gamma^{Z}(A)$, and hence, $\Gamma^{Z}(A) = \{A\cap C: C\in \Gamma^{Z}(P)\}$ holds.

$(5)\Rightarrow(3):$ Clearly, $D^{\delta}\subseteq D^{\delta}\mid_{x}$. Conversely, assume $m\in D^{\delta}\mid_{x}$, $n\in D^{u}$. Then $x, n\in D^{u}$. Since $D^{u}$ is filtered, there is a $p\in D^{u}$ such that $p\leq x,n$. This implies $p\in D^{u}\mid_{x}$, so $m\leq p$. It follows that $m\leq n$ by $p\leq n$. Thus $D^{\delta}\mid_{x}\subseteq D^{\delta}$.
\end{proof}

\begin{example}
The condition (5) in the above lemma is not equivalent to others. Let $\mathbb{N}$ be the set of natural numbers and $P = \mathbb{N}\cup(\mathbb{N}^{\partial}\dot{\cup}\mathbb{N}^{\partial}) $ with the partial order defined by $x\leq y$ iff $x\leq y$ in $\mathbb{N}$ or $x\leq y$ in $\mathbb{N}^{\partial}$ or $x\in \mathbb{N}$ and $y\in \mathbb{N}^{\partial}\dot{\cup}\mathbb{N}^{\partial} $ (see Fig. 1 for a better understanding).  One can easily sees that for any $-n\in \mathbb{N}^{\partial}\dot{\cup}\mathbb{N}^{\partial}$ and $D = \mathbb{N}\in \mathcal{D}({\downarrow}\{-n\})$, $D^{\delta}|_{-n} = D^{\delta}$, but $D^{u}$ is not filtered.

\begin{figure}[t]
	\centering
	\begin{tikzpicture}[scale=1.0]
		%\path (0,-1.2)  node[left] {} coordinate (a);
		\path (0,0)   node[left] {1} coordinate (b);
		\fill (b) circle (1pt);
		\path (0,1)   node[left] {2} coordinate (c);
		\fill (c) circle (1pt);
		\path (0,2)   node[left] {3} coordinate (d);
		\fill (d) circle (1pt);
		\path (0,2.8)  node[left] {} coordinate (j);
		\path (-1.5, 3.4)  node[left] {} coordinate (k);
		\path (1.5, 3.4)  node[left] {} coordinate (l);
		\path (-1.5,4.2)  node[left] {$-2$} coordinate (e);
		\fill (e) circle (1pt);
		\path (1.5,4.2)  node[right] {$-2$} coordinate (f);
		\fill (f) circle (1pt);
		\path (-1.5,5.2)  node[left] {$-1$} coordinate (g);
		\fill (g) circle (1pt);
		\path (1.5,5.2)  node[right] {$-1$} coordinate (h);
		\fill (h) circle (1pt);
		\draw   (b) -- (c);
		\draw  (c) -- (d) ;
		\draw      (e) -- (g)    (f) -- (h);
		\draw[densely dashed]   (d)--(j);
		\draw[densely dashed] (e)--(k);
		\draw[densely dashed]	(f)--(l)   ;
		
		\draw[densely dashed] (1,-0.4) arc (0:180:1 and 3.5);
	\end{tikzpicture}
	\caption{}
	%\captionsetup{labelsep=period}
\end{figure}
\end{example}

\begin{corollary}
Every zcpo $P$ has a lower hereditary $Z$-Scott topology.
\end{corollary}
\begin{proof}
Since $\sup D$ exists for any $D\in Z(P)$, $D^{u}$ is filtered.
\end{proof}

\begin{definition}
A poset $P$ is called locally weakly  meet $s_{Z}$-continuous if ${\downarrow}x$ as a subposet of $P$ is weakly  meet $s_{Z}$-continuous for each $x\in P$.
\end{definition}

\begin{lemma}\label{2l1}
Let $P$ be a poset with a lower hereditary $Z$-Scott topology and $A\in \Gamma^{Z}(P)$. Then for any $E\subseteq A$, we have $cl_{\sigma^{Z}(P)}(E) = cl_{\sigma^{Z}(A)}(E)$.
\end{lemma}
\begin{proof}
Straightforward.
\end{proof}

\begin{theorem}
Let $P$ be a poset with a lower hereditary $Z$-Scott topology and $Z$ be subset hereditary. Then $P$ is weakly meet $s_{Z}$-continuous if and only if $P$ is locally weakly meet $s_{Z}$-continuous.
\end{theorem}
\begin{proof}
$(\Rightarrow)$: For any $x\in P$, let $D\in Z(\da x)$, $y\in D^{\delta}\mid_{x}$. Then $D\in Z(P)$. Since the $Z$-Scott topology of $P$ is lower hereditary, by Lemma \ref{lh}, we have $D^{\delta}\mid_{x} = D^{\delta}$. Thus $y\in D^{\delta}$. It follows that $y\in cl_{\sigma^{Z}(P)}(\da y\cap \da D)$ by the weakly meet $s_{Z}$-continuity of $P$. Therefore, $y\in cl_{\sigma^{Z}(\da x)}(\da y\cap \da D)$ by Lemma \ref{2l1}.

$(\Leftarrow)$: Suppose $D\in Z(P)$, $y\in D^{\delta}$. For any $m\in D^{u}$, we have $D\subseteq \da m$ and $D\in Z(\da m)$ as $Z$ is subset hereditary. Since $D^{\delta} = D^{\delta}\mid_{m}$, we have $y\in D^{\delta}\mid_{m}$, which implies that $y\in cl_{\sigma^{Z}(\da m)}(\da y\cap \da D)$. So $y\in cl_{\sigma^{Z}(P)}(\da y\cap \da D)$ by Lemma \ref{2l1} again.
\end{proof}

\begin{proposition}\label{2p1}
Let $P$ be a poset with a lower hereditary $Z$-Scott topology. If $P$ is weak $s_{Z}$-continuous and for any $x\in P$, $y\in \da x$, $\ddo_{Z}^{x}y\in Z(\da x)$, then $\da x$ is $s_{Z}$-continuous.
\end{proposition}
\begin{proof}
We need to prove that $y\in (\ddo_{Z}^{x}y)^{\delta}\mid_{x}$. Since $P$ is weak $s_{Z}$-continuous, $y\in (\ddo_{Z}y)^{\delta}$. Assume $m\ll_{Z}y$, $D\in Z(\da x)$ with $y\in D^{\delta}\mid_{x}$. Then $y\in D^{\delta}$ since the $Z$-Scott topology on $P$ is lower hereditary. So $m\in \da D$ by $m\ll_{Z}y$, which implies that $m\ll_{Z}y$ in $\da x$. Therefore, $\ddo_{Z}y\subseteq \ddo_{Z}^{x}y$. Hence, $y\in (\ddo_{Z}y)^{\delta}\subseteq (\ddo_{Z}^{x}y)^{\delta} = (\ddo_{Z}^{x}y)^{\delta}\mid_{x}$, where the last equality holds as $\ddo_{Z}^{x}y\in Z(\da x)$. Moreover, $\ddo_{Z}^{x}y\in Z(\da x)$ implies $\ddo_{Z}^{x}y\in I_{Z}(\da x)$, so $\da x$ is $s_{Z}$-continuous.
\end{proof}

\begin{proposition}\label{2p2}
Let $P$ be a poset with a lower hereditary $Z$-Scott topology and $Z$ be subset hereditary. If for any $x\in P$, $\da x$ is $s_{Z}$-continuous and $\ddo_{Z}x\in Z(P)$, then $P$ is $s_{Z}$-continuous.
\end{proposition}
\begin{proof}
We only need to prove that $x\in (\ddo_{Z}x)^{\delta}$. By assumption, $\da x$ is $s_{Z}$-continuous, we have $x\in (\ddo_{Z}^{x}x)^{\delta}\mid_{x}$. Now we show that $\ddo_{Z}^{x}x\subseteq \ddo_{Z}x$. Let $m\ll_{Z}x$ in $\da x$ and $D\in Z(P)$ with $x\in D^{\delta}$. We can find that $D\subseteq \da y$ and $D\in Z(\da y)$ for each $y\in D^{u}$. Claim that $\ddo_{Z}^{x}x\subseteq \ddo_{Z}^{y}x$. Assume $a\in \ddo_{Z}^{x}x$. Since $\da y$ is $s_{Z}$-continuous and $x\in \da y$, we have $\ddo_{Z}^{y}x\in Z(\da y)$ and $x\in (\ddo_{Z}^{y}x)^{\delta}\mid_{y}$. It follows that $x\in (\ddo_{Z}^{y}x)^{\delta}\mid_{x}$ as $\da x\subseteq \da y$. Moreover, $\ddo_{Z}^{y}x\in Z(\da x)$ as $\ddo_{Z}^{y}x\subseteq \da x$ and $Z$ is subset hereditary. This implies that $a\in \ddo_{Z}^{y}x$. So $\ddo_{Z}^{x}x\subseteq \ddo_{Z}^{y}x$ holds. Thus $m\ll_{Z}x$ in $\da y$. As $x\in D^{\delta}$ implies that $x\in D^{\delta}\mid_{y}$, we have $m\in \da D$. Hence, $m\ll_{Z}x$. Then $\ddo_{Z}^{x}x\subseteq \ddo_{Z}x$. It is self-evident that $x\in (\ddo_{Z}^{x}x)^{\delta}\mid_{x}\subseteq (\ddo_{Z}x)^{\delta}\mid_{x} = (\ddo_{Z}x)^{\delta}$, where the last equality holds as $\ddo_{Z}x\in Z(\da x)$. So $P$ is $s_{Z}$-continuous.
\end{proof}

\begin{theorem}
Let $P$ be a poset with a lower hereditary $Z$-Scott topology and $Z$ be subset hereditary. Then the following conditions are equivalent:
\begin{enumerate}
\item[$\mathrm{(1)}$] $P$ is $s_{Z}$-continuous and $\ddo_{Z}^{x}y\in Z(\da x)$ for any $x\in P$ and $y\in \da x$;
\item[$\mathrm{(2)}$] $\da x$ is $s_{Z}$-continuous and $\ddo_{Z}x\in Z(P)$ for any $x\in P$.
\end{enumerate}
\end{theorem}
\begin{proof}
Straightforward by Proposition \ref{2p1} and \ref{2p2}.
\end{proof}

\section{A monad on POSET$_{\delta}$}

In this part, $\mathbf{POSET_{\delta}}$ denotes the category whose objects are all posets and morphisms are $\sigma^{Z}$-continuous mappings. We will give a monad on $\mathbf{POSET_{\delta}}$ and characterize its Eilenberg-Moore algebras.

\begin{definition}
Let $P$ be a poset and $x, y\in P$.
\begin{enumerate}
\item[$\mathrm{(1)}$] $x$ is called $Z$-$beneath$ $y$, denoted by $x\prec_{Z}y$, if for any $A\in \Gamma^{Z}(P)$ with $y\in A^{\delta}$, $x\in A$.
\item[$\mathrm{(2)}$] $P$ is said to be $\delta_{Z}$-$continuous$ if for all $a\in P$, $a\in \{m\in P: m\prec_{Z} a\}^{\delta}$.
\end{enumerate}
\end{definition}
Notice that the set $\{m\in P: m\prec_{Z} a\}\in \Gamma^{Z}(P)$ automatically. There are some common properties about the relation $\prec_{Z}$ being similar to the $\ll$.
\begin{proposition}
Let $P$ be a poset and $x, y, m, n\in P$.
\begin{enumerate}
 \item[$\mathrm{(1)}$] $x\prec_{Z} y$ implies $x\leq y$;
 \item[$\mathrm{(2)}$] $m\leq x\prec_{Z} y\leq n$ implies $m\prec_{Z} n$;
 \item[$\mathrm{(3)}$] if $P$ has a bottom $0$, then $0\prec_{Z} x$ always holds.
\end{enumerate}
\end{proposition}

\begin{proposition}\label{ccc}
Let $P$ be a poset and $\mathcal{C}\in\Gamma^{Z}(\Gamma^{Z}(P))$. Then the supremum of $\mathcal{C}$ in $\Gamma^{Z}(P)$ exists and is exactly $\bigcup \mathcal{C}$.
\end{proposition}
\begin{proof}
Clearly, it is enough to show that $\bigcup\mathcal{C}\in \Gamma^{Z}(P)$. For any $D\in Z(P)$ with $D\subseteq \bigcup\mathcal{C}$, there is $C_{d}\in \mathcal{C}$ for each $d\in D$ such that $d\in C_{d}$. Then we have ${\downarrow} d\subseteq C_{d}$ and $\{{\downarrow} d: d\in D\}\subseteq \mathcal{C}$ as $\mathcal{C}$ is a lower set. Since the monotonicity of the mapping $f: P\rightarrow \Gamma^{Z}(P)$ defined by $f(p) = {\downarrow} p$ implies that $\{{\downarrow} d: d\in D\}\in Z(\Gamma^{Z}(P))$, $\{{\downarrow} d: d\in D\}^{\delta}\subseteq \mathcal{C}$. Now consider each $a\in D^{\delta}$, we have ${\downarrow} a\in \{{\downarrow} d: d\in D\}^{\delta}$, which means ${\downarrow} a\in \mathcal{C}$. Thus $a\in \bigcup \mathcal{C}$, and $D^{\delta}\subseteq \bigcup \mathcal{C}$ holds.
\end{proof}

\begin{definition}
Let $P$ be a poset.
\begin{enumerate}
  \item[$\mathrm{(1)}$] An element $x$ of $P$ is called $Z$-$compact$ if $x\prec_{Z} x$. We use $k_{Z}(P)$ to denote the set of all $Z$-compact elements of $P$.
  \item[$\mathrm{(2)}$] $P$ is called $\delta_{Z}$-$prealgebraic$ if for each $x\in P$, $x\in \{y\in k_{Z}(P): y\leq x \}^{\delta}$.
\end{enumerate}
\end{definition}

Notably, we call a $\delta_{Z}$-$prealgebraic$ complete lattice a $\delta_{Z}$-$prealgebraic$ lattice for short. Obviously, $\Gamma^{Z}(P)$ is a $\delta_{Z}$-$prealgebraic$ lattice for any poset $P$.

\begin{lemma}\label{conn}
Let $(g, d)$ be a Galois connection between two posets $S$ and $T$, where $g: S\rightarrow T$, $d: T\rightarrow S$. Then $d$ preserves cuts of any subset of $T$, that is, $d(A^{\delta})\subseteq d(A)^{\delta}$ for any $A\subseteq T$.
\end{lemma}
\begin{proof}
It suffices to show that $d(x)\in d(A)^{\delta}$ for any $x\in A^{\delta}$. Let $y$ be an upper bound of $d(A)$. Then for each $a\in A$, we have $d(a)\leq y$, and so $a\leq g(y)$. It follows that $A\subseteq {\downarrow}g(y)$. Thus $A^{\delta}\subseteq {\downarrow}g(y)$, which implies $x\leq g(y)$, so $d(x)\leq y$. Hence, $d(x)\in d(A)^{\delta}$.
\end{proof}

\begin{lemma}\label{close}
Let $(g, d)$ be a Galois connection between two posets $S$ and $T$, where $g: S\rightarrow T$, $d: T\rightarrow S$. Then for any $C\in \Gamma^{Z}(S)$, ${\downarrow}g(C)\in \Gamma^{Z}(T)$.
\end{lemma}
\begin{proof}
Let $E$ be a $Z$-set of $T$ with $E\subseteq {\downarrow}g(C)$.  Then for any $e\in E$, there is a  $c_{e}\in C$ such that $e\leq g(c_{e})$, this means $d(e)\leq c_{e}$. Thus $d(E)\subseteq C$ and $d(E)^{\delta}\subseteq C$ since $C\in \Gamma^{Z}(S)$. The conclusion of Lemma \ref{conn} indicates that $d(E^{\delta})\subseteq C$. Therefore, $E^{\delta}\subseteq d^{-1}(C) = {\downarrow}g(C)$.
\end{proof}

\begin{lemma}\label{dgg}
Let $(g, d)$ be a Galois connection between two posets $S$ and $T$, where $g: S\rightarrow T$, $d: T\rightarrow S$. Consider the following two conditions:
\begin{enumerate}
\item[$\mathrm{(1)}$] For any $A\in \Gamma^{Z}(S)$, $g(A^{\delta})\subseteq g(A)^{\delta}$.
\item[$\mathrm{(2)}$] $d$ preserves the relation $\prec_{Z}$.
\end{enumerate}
Then $(1)\Rightarrow (2)$; if $T$ is $\delta_{Z}$-continuous, then $(2)\Rightarrow (1)$.
\end{lemma}
\begin{proof}
$(1)\Rightarrow (2)$: We need to show that $d(x)\prec_{Z}d(y)$ for any $x\prec_{Z} y$ in $T$. Let $A\in \Gamma^{Z}(S)$ with $d(y)\in A^{\delta}$. Then $y\leq g(m)$ for some $m\in A^{\delta}$. By the condition (1), we have  $g(m)\in g(A)^{\delta}$ and hence, $y\in g(A)^{\delta} = ({\downarrow g(A)})^{\delta}$. Lemma \ref{close} indicates that ${\downarrow}g(A)\in \Gamma^{Z}(T)$, then $x\in {\downarrow}g(A)$ as $x\prec_{Z}y$. Thus there is an $a\in A$ such that $x\leq g(a)$, which implies $d(x)\leq a$. It follows that $d(x)\in A$. Therefore, $d$ preserves the relation $\prec_{Z}$.

$(2)\Rightarrow (1)$: By the $\delta_{Z}$-continuity of $T$, we know $g(x)\in \{y\in T: y\prec_{Z}g(x)\}^{\delta}$ for  each $x\in A^{\delta}$. Thus in order to show $g(A^{\delta})\subseteq g(A)^{\delta}$ for any $A\in \Gamma^{Z}(S)$, it suffices to prove that for each $x\in A^{\delta}$, $\{y\in T: y\prec_{Z}g(x)\}^{\delta}\subseteq g(A)^{\delta}$. For each $y\prec_{Z}g(x)$, since $d(y)\prec_{Z}d(g(x))\leq x$, we have $d(y)\prec_{Z} x$. Then $d(y)\in A$ because $x\in A^{\delta}$ and $A\in \Gamma^{Z}(S)$, which implies $y\in {\downarrow}g(A)$. Thus $\{y\in T: y\prec_{Z}g(x)\}^{\delta}\subseteq ({\downarrow} g(A))^{\delta} = g(A)^{\delta}$.
\end{proof}

\begin{lemma}\label{zzz}
If $L$ is a $zcpo$, then $k_{Z}(L)$ is also a $zcpo$.
\end{lemma}
\begin{proof}
We just need to prove that $\sup D\in k_{Z}(L)$ for any $D \in Z(k_{Z}(L))$. Let $A\in \Gamma^{Z}(L)$ with $\sup D\in A^{\delta}$. Then $D\subseteq A^{\delta}$, and so $D\subseteq A$ by $D\subseteq k_{Z}(L)$. Thus ${\downarrow}\sup D = D^{\delta}\subseteq A$, this means $\sup D\in A$. Hence, $\sup D\in k_{Z}(L)$ and $\sup_{k_{Z}(L)}D = \sup D\in k_{Z}(L)$. It follows that $k_{Z}(L)$ is a $zcpo$.
\end{proof}

The above lemma ensures that $D^{\delta}\mid_{k_{Z}(L)} = D^{\delta}\cap k_{Z}(L)$ hold. There is an example illustrating that $D^{\delta}\mid_{k_{Z}(L)} = D^{\delta}\cap k_{Z}(L)$ doesn't hold  if $L$ is not a $zcpo$.

\begin{example}
Let $P$ be the poset consist of all natural numbers $\mathbb{N}$ and $\{a, b, c, d, \top\}$. $\top$ is the greatest element of $P$ and $\{a, b, c, d\}\subseteq \mathbb{N}^{u}$, $c\in \{a, b\}^{u}$. Now consider $Z = D$, where $D(P)$ is the family of all directed subsets. It is easy to verify that $k_{D}(P) = \mathbb{N}\cup \{d\}$. For $\mathbb{N}\in D(P)$, $\mathbb{N}^{\delta}\mid_{k_{D}(P)} = \mathbb{N}\cup \{d\}$, however, $\mathbb{N}^{\delta}\cap k_{D}(P) = \mathbb{N}$ since $\mathbb{N}^{\delta} = \mathbb{N}$.
\end{example}

 We denote by $\mathbf{\delta_{Z}PALG}$ the category which has all $\delta_{Z}$-prealgebraic lattices as objects and  maps that have an upper adjoint and preserve the relation $\prec_{Z}$ as morphisms. Next, we will investigate the relation between the categories $\mathbf{POSET_{\delta}}$ and $\mathbf{\delta_{Z}PALG}$.

\begin{theorem}
Let $K_{Z}$ and $\Gamma^{Z}$ be two functors between $\mathbf{\delta_{Z}PALG}$ and $\mathbf{POSET_{\delta}}$. Here $K_{Z}$ is defined by associating a $\delta_{Z}$-prealgebraic lattice with the poset $k_{Z}(L)$ and a morphism $f: L \rightarrow M$ in $\mathbf{\delta_{Z}PALG}$ with the map $K_{Z}(f): k_{Z}(L)\rightarrow k_{Z}(M)$ defined by
\begin{center}
 $\forall x\in k_{Z}(L)$, $K_{Z}(f)(x) = f(x)$;
\end{center}
$\Gamma^{Z}$ is defined by assigning a poset $P$ to the $\delta_{Z}$-prealgebraic lattice $\Gamma^{Z}(L)$ and the $\sigma^{Z}$-continuous mapping $g: P\rightarrow Q$ to $\Gamma^{Z}(g): \Gamma^{Z}(P)\rightarrow \Gamma^{Z}(Q)$ defined as follows:
\begin{center}
 $\forall A\in \Gamma^{Z}(P)$, $\Gamma^{Z}(g)(A) = cl_{\sigma^{Z}(Q)}(g(A))$.
\end{center}
Then $\Gamma^{Z}$ is left adjoint to $K_{Z}$ with unit $\eta_{P}$ and counit $\epsilon_{P}$ given by
\begin{center}
 $\eta_{P}: P\rightarrow K_{Z}\Gamma^{Z}(P): p\mapsto {\downarrow}p$, $\forall p\in P$, and
\end{center}
\begin{center}
$\epsilon_{L}: \Gamma^{Z}K_{Z}(L)\rightarrow L: E\mapsto \sup E$, $\forall E\in \Gamma^{Z}K_{Z}(L)$,
\end{center}
respectively.
\end{theorem}
\begin{proof}
Step 1:  Verify that functors $K_{Z}$ and $\Gamma^{Z}$ are well-defined by showing that $K_{Z}(f)$ and $\Gamma^{Z}(g)$ are morphisms in $\mathbf{POSET_{\delta}}$, $\mathbf{\delta_{Z}PALG}$, respectively. We claim that $K_{Z}(f)$ is $\sigma^{Z}$-continuous, that is, $K_{Z}(f)(D^{\delta}\mid_{k_{Z}(L)} )\subseteq (K_{Z}(f)(D))^{\delta}\mid _{k_{Z}(M)}$ for any $D\in Z(k_{Z}(L))$. Since $L$ is a complete lattice, by Lemma \ref{zzz}, $k_{Z}(L)$ is a $zcpo$. Thus we only need to prove that $K_{Z}(f)({\downarrow}_{k_{Z}(L)} \sup_{k_{Z}(L)}D)\subseteq {\downarrow}_{k_{Z}(M)}\sup_{k_{Z}(M)}K_{Z}(f)(D)$. More precisely, to show $K_{Z}(f)({\downarrow}\sup D\cap k_{Z}(L))\subseteq {\downarrow}\sup K_{Z}(f)(D)\cap k_{Z}(M)$. From the fact that $f$ has an upper adjoint, we know $f(\sup D) = \sup f(D)$ holds. So it is easy to see that
\begin{align*}
  K_{Z}(f)({\downarrow}\sup D\cap k_{Z}(L)) &\subseteq K_{Z}(f)({\downarrow}\sup D)\cap k_{Z}(M) \\
  & \subseteq {\downarrow}f(\sup D)\cap k_{Z}(M) \\ &= {\downarrow}\sup f(D)\cap k_{Z}(M) \\ & = {\downarrow}\sup K_{Z}(f)(D)\cap k_{Z}(M).
\end{align*}
Hence, $K_{Z}(f)$ is $\sigma^{Z}$-continuous.

We proceed to show $\Gamma^{Z}(g)$ has an upper adjoint and preserves the relation $\prec_{Z}$. It is obvious that $\Gamma^{Z}(g)$ preserves arbitrary sups in $\Gamma^{Z}(P)$, by Corollary $\mathrm{O}${-}3.5 in \cite{clad}, $\Gamma^{Z}(g)$ has an upper adjoint. Moreover, the upper adjoint is given by
\begin{center}
$h: \Gamma^{Z}(Q)\rightarrow \Gamma^{Z}(P): C\mapsto  g^{-1}(C)$.
\end{center}
By Proposition \ref{ccc}, we know for any $\mathcal{C}\in \Gamma^{Z}(\Gamma^{Z}(Q)), \sup \mathcal{C} = \bigcup \mathcal{C}$. It follows that
\begin{center}
$h(\mathcal{C}^{\delta}) = g^{-1}({\downarrow}\sup \mathcal{C}) = {\downarrow}\{g^{-1}(\bigcup \mathcal{C})\} = {\downarrow}\{\bigcup g^{-1}(\mathcal{C})\} = {\downarrow}\sup h(\mathcal{C}) = h(\mathcal{C})^{\delta}$.
\end{center}
Therefore, from the conclusion of Lemma \ref{dgg}, we get that $\Gamma^{Z}(g)$ preserves $\prec_{Z}$. So $\Gamma^{Z}(g)$ is a morphism in $\mathbf{\delta_{Z}PALG}$.

Step 2: To show $\Gamma^{Z}$ is left adjoint to $K_{Z}$ in detail. Obviously, $\eta_{P}$ is $\sigma^{Z}$-continuous, that is, a morphism in $\mathbf{POSET_{\delta}}$. Now let $L$ be a $\delta_{Z}$-prealgebraic lattice and $f: P\rightarrow K_{Z}(L)$ $\sigma^{Z}$-continuous. We define
\begin{center}
$\bar{f}: \Gamma^{Z}(P)\rightarrow L: A\mapsto \sup f(A)$.
\end{center}
It is easy to find that $K_{Z}(\bar{f})\circ \eta = f$. Thus for the remainder, we need to prove that $\bar{f}$ is a unique morphism in $\mathbf{\delta_{Z}PALG}$ such that $K_{Z}(\bar{f})\circ \eta = f$. Note that $\bar{f}$ preserves arbitrary sups in $\Gamma^{Z}(P)$ and $\Gamma^{Z}(P)$ is a complete lattice, so $\bar{f}$ has an upper adjoint, denoted by $f^{*}$. More specifically, for any $m\in L$,

\begin{align*}
f^{*}(m) &= \sup\bar{f}^{-1}({\downarrow}m) \\ &= \sup\{C\in \Gamma^{Z}(P): \bar{f}(C)\leq m\} \\ &= \sup\{{\downarrow}t\in \Gamma^{Z}(P): \bar{f}({\downarrow}t)\leq m\} \\ &= \sup\{{\downarrow}t\in \Gamma^{Z}(P): f(t)\leq m\}.
\end{align*}
Then again by Lemma \ref{dgg}, we check that $f^{*}(B^{\delta})\subseteq f^{*}(B)^{\delta}$ for any $B\in \Gamma^{Z}(L)$ to affirm $\bar{f}$ preserves $\prec_{Z}$. Since $L$ and $\Gamma^{Z}(P)$ are complete lattices, we only need to prove $f^{*}(\sup B)\leq \sup f^{*}(B)$. To this end, consider each ${\downarrow}x\in \Gamma^{Z}(P)$ which satisfies $f(x)\leq \sup B$, that is, $f(x)\in B^{\delta}$. Then $f(x)\in B$ as $f(x)\in k_{Z}(L)$. It follows that ${\downarrow}x\subseteq f^{*}(f(x))$, in addition, ${\downarrow}x\subseteq\sup f^{*}(B)$. Thus $f^{*}(\sup B)\leq \sup f^{*}(B)$ holds, hence, $f^{*}(B^{\delta})\subseteq f^{*}(B)^{\delta}$. Besides, clearly, $\bar{f}$ is unique. Therefore, we can conclude that $\Gamma^{Z}$ is left adjoint to $K_{Z}$.
\end{proof}

Next, we will give a monad on $\mathbf{POSET_{\delta}}$. Before this, let us recall the following conclusion:
\begin{proposition}\cite{ct}
Let $U: \mathcal{B}\rightarrow \mathcal{A}$ and $F: \mathcal{A}\rightarrow \mathcal{B}$ be functors such that $F$ is left adjoint to $U$ with $\eta: id\rightarrow UF$ and $\epsilon: FU\rightarrow id$ the unit and counit, respectively. Then $(UF, \eta, U\epsilon F)$ is a monad on $\mathcal{A}$.
\end{proposition}
Now, by combining the above two conclusions, and $K_{Z}\Gamma^{Z}$ is written as $\delta$, we obtain the following.
\begin{theorem}
The endofunctor $\delta$  together with two natural transformation $\eta: id\rightarrow \delta$ and $\mu = \Gamma^{Z}\epsilon K_{Z}: \delta^{2}\rightarrow\delta$ is a monad on the category $\mathbf{POSET_{\delta}}$. More precisely, for each $P\in \mathbf{POSET_{\delta}}$, $\eta_{P}: P\rightarrow \delta(P)$ and $\mu_{P}: \delta^{2}(P)\rightarrow \delta(P)$ are defined as:
\begin{center}
$\forall p\in P, \eta(p) = {\downarrow}p$,
\end{center}
\begin{center}
$\forall \mathcal{A}\in \delta^{2}(P), \mu(\mathcal{A}) = \sup \mathcal{A}$,
\end{center}
respectively.
\end{theorem}

Recall that an $Eilenberg$-$Moore$ algebra for a monad $(T, \eta, \mu)$ on a category $\mathcal{C}$ is a pair $(C, \xi)$, where $\xi: TC\rightarrow C$ is a morphism in $\mathcal{C}$ called a structure map which satisfies $\xi\circ\eta_{C} = id_{C}$ and $\xi\circ\mu_{C} = \xi T\xi$. In addition, we call a poset $P$ $\delta cpo$ if for any $A\in \delta(P), \sup A$ exists.

\begin{theorem}
There exists a structure map $\xi: \delta(P)\rightarrow P$ in $\mathbf{POSET_{\delta}}$ such that $(P, \xi)$ is an $Eilenberg$-$Moore$ algebra of the monad $(\delta, \eta, \mu)$ if and only if $P$ is a $\delta cpo$.
\end{theorem}
\begin{proof}
$(\Rightarrow)$: We claim that $\sup A = \xi(A)$ exists for any $A\in \delta(P)$. Consider each $a\in A$, from the facts that $\eta(a) = {\downarrow}a\subseteq A$ and $\xi$ is order-preserving, we have $\xi(\eta(a)) \leq \xi(A)$. This implies $a\leq \xi(A)$ since $\xi\circ\eta_{P} = id_{P}$. Thus $\xi(A)$ is an upper bound of $A$. Assume that $m$ is another upper bound of $A$, which means $A\subseteq {\downarrow}m = \eta(m)$. So $\xi(A)\leq\xi(\eta(m))$ by the monotonicity of $\xi$ again. It follows that $\xi(A)\leq m$, and hence $\xi(A) = \sup A$.

$(\Leftarrow)$: Since $P$ is a $\delta cpo$, we can define $\xi: \delta(P)\rightarrow P$ by $C\mapsto \sup C$. One can easily verify that $\xi\circ\eta_{P} = id_{P}$ and $\xi\circ\mu_{P} = \xi\delta\xi$. For the remaining part, what we need to prove is that $\xi$ is a morphism in $\mathbf{POSET_{\delta}}$, that is, $\xi$ is $\sigma^{Z}$-continuous. To this end, let $\mathcal{A}$ be a $Z$-set of $\delta(P)$. Since for every upper bound $y$ of $\xi(\mathcal{A})$, we have $\xi(A) = \sup A\leq y$ for each $A\in \mathcal{A}$, which implies $A\subseteq {\downarrow}y$ and so ${\downarrow}y$ is an upper bound of $\mathcal{A}$. Thus for every $B\in \mathcal{A}^{\delta}$, $B\subseteq{\downarrow}y$, which means $\xi(B) = \sup B\leq y$. It follows that $\xi(\mathcal{A}^{\delta}) = \{\xi(B): B\in \mathcal{A}^{\delta}\}\subseteq \xi(\mathcal{A})^{\delta}$. Hence, $\xi$ is $\sigma^{Z}$-continuous.
\end{proof}

Refer to \cite{andrea}, if $(T, \eta, \mu)$ is a monad on the category $\mathcal{C}$, the category $\mathbf{T}${-}$\mathbf{ALG}$ consists of all $T$-algebras as objects and morphisms of $T$-algebras as morphisms. Here a morphism of $T$-algebras between $(C, \xi)$ and $(C', \xi')$ in $\mathbf{T}${-}$\mathbf{ALG}$ is a morphism $f: C\rightarrow C'$ in $\mathcal{C}$ which satisfies $f\circ\xi = \xi'Tf$.

Combining with the characterization of $\delta$-algebras, we can deduce that $f: (P, \alpha)\rightarrow (Q, \beta)$ is a $\delta$-algebra morphism if and only if $f(\sup A) = \sup f(A)$ for every $A\in \delta(P)$.

%\section{references}%\label{references}

\bibliographystyle{./entics}

\begin{thebibliography}{10}\label{bibliography}



\bibitem{czp}
Bandelt, H.J., M. Ern$\acute{e}$, \emph{The category of Z-continuous posets},
 Journal of Pure and Applied Algebra, \textbf{30} (1983): 219-226.
\newline\url{https://doi.org/10.1016/0022-4049(83)90057-9}

\bibitem{zcp}
Baranga, A.,  \emph{$Z$-continuous posets},
 Discrete Mathematics, \textbf{152} (1996): 33-45.
\newline \url{https://doi.org/10.1016/0012-365X(94)00307-5}


\bibitem{ct}
Barr, M. and C. Wells, ``Category Theory Lecture Notes for ESSLLI", Lecture Notes for ESSLLI, 1999. Available online at:
\newline \url{https://fldit-www.cs.tu-dortmund.de/~peter/barrwells.pdf}%https://doi.org/10.1017/CBO9781139524438}

\bibitem{}
 Davey, B.A. and H.A. Priestley, ``Introduction to lattices and order", Cambridge University Press, 2002.
ISBN: 9780521784511



\bibitem{erne}
Ern$\acute{e}$, M., \emph{Scott convergence and Scott topology in partially ordered sets II},
Continuous lattices. Springer, Berlin, Heidelberg, (1981): 61-96.
\newline \url{https://doi.org/10.1007/BFb0089904}

\bibitem{clad}
Gierz, G., K.~H. Hofmann, K.~Keimel, J.~D. Lawson, M.~Mislove, and D.~S. Scott, ``Continuous Lattices and Domains", Volume~\textbf{93} of Encyclopedia of Mathematics and its Applications, Cambridge University Press, 2003.
\newline \url{https://doi.org/10.1017/CBO9780511542725}

\bibitem{nht2}
Goubault-Larrecq,  J.,``Non-Hausdorff Topology and Domain Theory", Volume~\textbf{22} of New Mathematical Monographs, Cambridge University Press, 2013.
\newline \url{https://doi.org/10.1017/CBO9781139524438}


\bibitem{zds}
Ho, W.K. and D.S. Zhao, \emph{Lattices of Scott-closed sets},
 Commentationes Mathematicae Universitatis Carolinae, \textbf{50} (2009): 297-314. Available online at:\newline
 \url {http://dml.cz/dmlcz/133435}

\bibitem{mp}
 Mao, X.X. and  L.S. Xu,
 \emph{Meet continuity properties of posets},
 Theoretical Computer Science \textbf{410} (2009): 4234-4240.
\newline \url{https://doi.org/10.1016/j.tcs.2009.06.017}


\bibitem{sz}
 Ruan, X.J. and X.Q. Xu,
 \emph{$s_{Z}$-Quasicontinuous posets and meet $s_{Z}$-continuous posets},
 Topology and Its Applications \textbf{230} (2017): 295-307.
\newline \url{https://doi.org/10.1016/j.topol.2017.08.045}


\bibitem{andrea}
Schalk, A., 
 \emph{Algebras for generalized power constructions},
 Technische Hochschule Darmstadt, Doctoral thesis, 1993. Available online at:
\newline \url{https://www.cs.man.ac.uk/~schalk/publ/diss.ps.gz}



\bibitem{ipic}
Wright, J.B., E.G. Wagner, J.W. Thatcher, \emph{A uniform approach to inductive posets and inductive closure},
 Theoretical Computer Science, \textbf{7} (1978): 57-77.
\newline \url{https://doi.org/10.1016/0304-3975(78)90040-3}



\bibitem{qzm}
Xu, X.Q.,  M.K. Luo and Y. Huang, \emph{Quasi Z-continuous domains and Z-meet continuous domains},
 Acta Mathematica Sinica. Chinese Series, \textbf{48} (2005): 221-234.
\newline \url{https://doi.org/10.3321/j.issn:0583-1431.2005.02.002} (Chinese) and \href{https://actamath.cjoe.ac.cn/Jwk_sxxb_cn/EN/10.12386/A20050026}{English listing}


\bibitem{s2}
 Zhang, W.F. and X.Q. Xu, \emph{$s_{2}$-Quasicontinuous posets},
 Theoretical Computer Science \textbf{574} (2015): 78-85.
\newline \url{https://doi.org/10.1016/j.tcs.2015.01.037}



\bibitem{cscp}
 Zhao, D.S., \emph{Closure spaces and completions of posets},
 Semigroup Forum \textbf{90} (2015): 545-555.
\newline \url{https://doi.org/10.1007/s00233-015-9692-6}
\end{thebibliography}

\end{document}